\documentclass[conference]{IEEEtran}
\IEEEoverridecommandlockouts
\usepackage{cite}
\usepackage{amsmath,amssymb,amsfonts}
\usepackage{graphicx}
\usepackage{textcomp}
\usepackage{xcolor}
\def\BibTeX{{\rm B\kern-.05em{\sc i\kern-.025em b}\kern-.08em
    T\kern-.1667em\lower.7ex\hbox{E}\kern-.125emX}}

\usepackage{bm}
\usepackage{algorithm}
\usepackage{algpseudocode}
\usepackage{enumitem}

\begin{document}

\title{A Novel Cryptography-Based Privacy-Preserving Decentralized Optimization Paradigm\\
}

\author{\IEEEauthorblockN{1\textsuperscript{st} Xiang Huo }
\IEEEauthorblockA{\textit{Department of Electrical and Computer Engineering} \\
\textit{University of Utah}\\
Salt Lake City, USA \\
xiang.huo@utah.edu}
\and
\IEEEauthorblockN{2\textsuperscript{nd} Mingxi Liu}
\IEEEauthorblockA{\textit{Department of Electrical and Computer Engineering} \\
\textit{University of Utah}\\
Salt Lake City, USA \\
mingxi.liu@utah.edu}
}

\maketitle

\begin{abstract}
Existing large-scale optimization schemes are challenged by both scalability and cyber-security. With the favorable scalability, adaptability, and flexibility, decentralized and distributed optimization paradigms are widely adopted in cyber-physical system applications. However, most existing approaches heavily rely on explicit information exchange between agents or between agents and the system operator, leading the entire framework prone to privacy risks. To tackle this issue, this paper synthesizes cryptography and decentralized optimization techniques to develop a novel privacy-preserving decentralized optimization paradigm. The proposed paradigm is generically applicable to strongly coupled convex optimization problems with nonseparable objective functions and linearly coupled constraints. The security and accuracy of the proposed paradigm are verified via numerical examples.

\end{abstract}


\section{Introduction}

Large-scale optimization problems broadly exist in industrial cyber-physical system (ICPS) applications, e.g., economic dispatch in power system \cite{mao2020privacy}, traffic congestion control in transportation system \cite{he2007towards}, distributed energy resource control in smart grid \cite{dobbe2019toward}, and resource allocation \cite{chen2017stochastic}. Owing to the scalability and flexibility, decentralized optimization has been extensively adopted to speed up large-scale optimization.

In most decentralized and distributed optimization approaches \cite{peng2014distributed,boyd2011distributed,rivera2016distributed}, agents and the system operator (SO) are mandated to exchange and expose private information, e.g., decision variables and the network topology information. This can lead to undesired privacy leakage, i.e., adversaries who wiretap the communication may steal private information of the agents or the SO \cite{zhang2018enabling}. Besides, agents may not be willing to share their personal information to any neighbors or the SO \cite{lu2018privacy}. Having these privacy issues identified, privacy-preserving strategies for decentralized and distributed optimization are drawing more attention. Differential privacy (DP) \cite{dwork2006calibrating}, which protects the participants' privacy by adding carefully-designed noises to the transmitted messages, has emerged to be a popular method \cite{han2016differentially,zhang2016dynamic}. An $\epsilon$-differential privacy-preserving algorithm was developed in \cite{han2016differentially} for distributed constrained optimization by adding perturbations to the public signal. To minimize the accuracy deterioration caused by perturbations, Zhang and Zhu \cite{zhang2016dynamic} analyzed the privacy-accuracy tradeoff and proposed an optimal dynamic DP mechanism with guidelines to choose privacy parameters. Unfortunately, the privacy-accuracy tradeoff fundamentally exists in DP based schemes \cite{nozari2016differentially}. Similar to DP, the emerging correlated noise based obfuscation strategy uses obfuscation to mask the true state value. Under this umbrella, a stochastic gradient descent based algorithm was develop in \cite{gade2018privacy} to protect the privacy via gradient obfuscation. Mak \emph{et al.} \cite{mak2020privacy} proposed a differentially private distributed algorithm with high fidelity and obfuscation quality based on the Alternating Direction Method of Multipliers (ADMM) algorithm. Nonetheless, obfuscation based methods invariably suffer from accuracy loss and therefore could lead to severe feasibility and optimality issues.

Privacy preservation can also be enabled by cryptography which obscures critical information in ciphertext. As decentralized algorithms normally involve operations on the transmitted messages, such operations on the ciphertexts are imperative in any cryptography-based decentralized algorithm. Among various cryptosystems, homomorphic encryption \cite{rivest1978data} allows arithmetic operations over ciphertext, leading to a result in an encrypted form that, when decrypted, matches the result of operations performed on the plaintext. Partially homomorphic encryption was incorporated into the projected subgradient algorithm \cite{zhang2018enabling} and ADMM \cite{zhang2018admm} to preserve the privacy of both intermediate states and gradients. However, neither \cite{zhang2018enabling} nor \cite{zhang2018admm} is capable of solving strongly coupled optimization problems with nonseparable objective functions and coupled constraints. As an improvement, Lu and Zhu \cite{lu2018privacy} designed a private key secure computation scheme which is applicable to general projected gradient-based algorithms. However, this scheme is not semantically secure. Besides, the proposed private key algorithm requires fully homomorphic encryption, and the public key encryption scheme is not implementable without sacrificing privacy. Moreover, several key assumptions of the security modeling in \cite{lu2018privacy}, e.g., agents cannot eavesdrop others' communication links and the SO can only launch temporarily independent attacks, have excluded some high-probability privacy events in real ICPS applications.

In this paper, we aim to improve the technical baseline of privacy-preserving decentralized algorithms through cryptography. Particularly, we consider a set of agents and a SO that securely compute a strongly coupled optimization problem. The contribution of this paper is four-fold: (1) We synthesize cryptography and decentralized optimization to develop a novel privacy-preserving scalable optimization paradigm for strongly coupled optimization problems. (2) An idea of coefficient assignment is proposed to preserve the general privacy of both the agents and the SO. Comparing with \cite{lu2018privacy}, the proposed method is not confined to only protecting the ``coefficient'' of the gradient. (3) The enhanced security that is achieved by eliminating assumptions on adversarial behaviors has not been tackled before. (4) Demands on the cryptosystem are relaxed -- only additively homomorphic encryption is required. Besides, both private key and public key based encryption schemes can be applied in the proposed privacy-preserving paradigm.

\section{Main Results}
\subsection{Problem Formulation}
This paper considers a family of strongly coupled optimization problems with nonseparable convex objective functions, globally coupled constraints, and local constraints. In specific, the optimization problem can be formulated as 
\begin{equation}
\begin{aligned}
& \underset{x}{\text{min}} & & \mathcal{F}(x) \\
& \: \text{s.t.} & & x_{i} \in \mathbb{X}_{i}, \quad \forall i=1,\cdots,n, \\
& & &  h(x) \leq 0,
\end{aligned}
\label{4}
\end{equation}
where $n$ agents are involved in this problem, $x_i$ denotes the decision variable of the $i$th agent, and $x=[x_1\cdots x_n]^{\mathsf{T}}$. The objective function considered here can be written as
\begin{align}
    {\mathcal{F}(x)} \triangleq \frac{1}{2} \|\sum_{i=1}^{n} A_{ui} x_i + c \|^2_2
    +\sum_{i=1}^{n}f_i(x_i),
    \label{1}
\end{align}
where the quadratic term denotes the coupled global objective and $f_i(x_i)$ is agent $i$'s local objective. The local constraint of agent $i$ is represented as a convex set $\mathbb{X}_i$. The global constraint is considered to be linearly coupled, which can be written as
\begin{equation}
    h(x) \triangleq \sum_{i=1}^{n} A_{gi} x_i + d \leq 0. 
    \label{2}
\end{equation}
The relaxed Lagrangian of \eqref{4} can be calculated as
\begin{equation} \label{5}
    \mathcal{L}(x,\lambda)=\mathcal{F}(x)+\lambda^{\mathsf{T}}h(x), 
\end{equation}
where $\lambda$ is the dual variable. Consequently, the subgradients of $\mathcal{L}(x,\lambda)$ \emph{w.r.t.} $x_i$ and $\lambda$ can be calculated as
\begin{subequations}
\begin{align}
    \nabla_{x_i} \mathcal{L}(x,\lambda) &{=} {A_{ui}^{\mathsf{T}}} (\sum_{i=1}^{n} A_{ui} x_i + c )
    {+}\nabla_{x_i} f_i(x_i) {+}  {A_{gi}}^\mathsf{T}\lambda, \label{6}\\
     \nabla_{\lambda} \mathcal{L}(x,\lambda) &{=} \sum_{i=1}^{n} A_{gi} x_i + d.    \label{7}
\end{align}
\end{subequations}

 These subgradients are critical for the updates in general primal-dual based algorithms \cite{koshal2011multiuser,liu2013decentralized}. However, because the coupling terms $\sum_{i=1}^{n} A_{ui} x_i$ in \eqref{6} and  $\sum_{i=1}^{n} A_{gi} x_i$ in \eqref{7} contain the aggregated information of all agents, intermediate decision variables of the agents must be transmitted to the SO for primal and dual updates, putting agents' privacy at risk. Besides, $c$ in \eqref{6} and $d$ in \eqref{7} could be private information held by either the agents or the SO depending on particular applications, leading to different situations in privacy preservation. This paper develops a privacy-preserving paradigm that protects the privacy of all the participants in solving problem \eqref{4} and enjoys the scalability of the decentralized optimization.

\subsection{Privacy-Preserving Paradigm Design}

The privacy preservation should be considered in two aspects. First, information exchanged between the agents and the SO must not reveal useful information to external eavesdroppers. Second, all participants should only have access to what they need to know. In solving \eqref{4}, private information includes the coefficients in the optimization problem and private decision variables of the agents. Specifically, $A_{u1},\ldots,A_{un}$, and $c$ in \eqref{1} are aggregated global objective coefficients and $A_{g1},\ldots,A_{gn}$, and $d$ in \eqref{2} are the aggregated global constraint coefficients. These global coefficients can be assigned to either the SO or the agents depending on the particular application. Note that the privacy situation varies depending on the role of the SO, i.e., a pure computation agency that has no knowledge about the optimization problem \cite{roy2012security} or a participant with its own private information. Therefore, based on whether the SO is a trustworthy third party, two cases exist:
\begin{itemize}
  \item {\textbf{Case 1}:} SO is trustworthy and has access to the global coefficients.
  \item {\textbf{Case 2}:} SO is blind from the optimization problem and does not hold any private information.
\end{itemize}

Due to the limit of space, we only focus on Case 1 in this paper. Case 2 will be discussed in a journal extension of this work. In Case 1, the SO possesses a set of private coefficients for solving the optimization problem. For example, for the load aggregation problem in power system, SO holds global coefficients $d$ and $c$ as its private information which could represent the load capacity and baseline load, respectively. Such information is considered proprietary and should be kept only to the SO to prevent malicious attacks on the distribution network.
Therefore, we consider the SO holding $c$, $d$, and all global coefficients $A_{ui},A_{gi}, \forall i=1,\cdots,n$ as its private information. The $i$th agent holds the global coefficients related to itself, i.e.,  $A_{ui}$ and $A_{gi}$. Besides,  the $i$th agent should also have a set of coefficients included in its local objective $f_i(x_i)$ denoted by $\mathcal{I}_{fi}$ and the private key $w$.

\begin{figure}[!t]
    \centering
    \includegraphics[width=0.41\textwidth,trim = 0mm 20mm 0mm 20mm, clip]{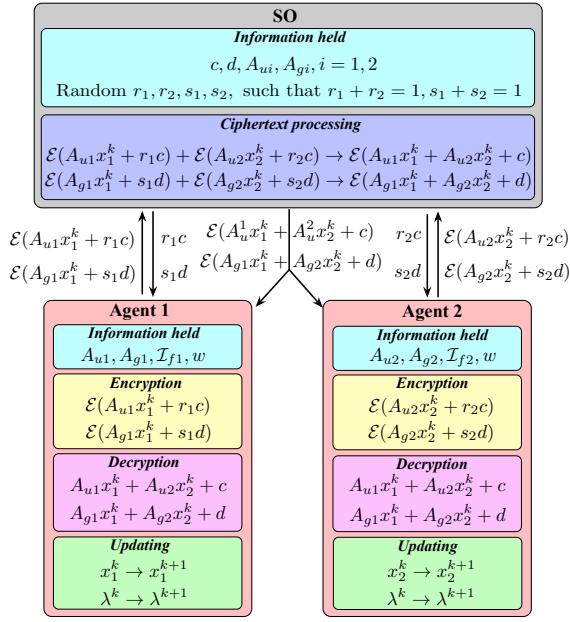}
    \caption{Cryptographic information flow between two agents and one SO.}
    \label{information_flow}
\end{figure}

We then propose a novel cryptography-based decentralized algorithm for problem \eqref{4}. Fig. \ref{information_flow} presents an illustrative example with two agents and one SO. During each iteration of any primal-dual based decentralized algorithm,
the SO randomly generates two sets of parameters, satisfying 
\begin{subequations}
\begin{align}
    r_1+r_2+\cdots+r_n &= 1, \label{8a}\\
    s_1+s_2+\cdots+s_n &= 1,
\end{align}
\end{subequations}
and sends both $r_ic$ and $s_id$ to the $i$th agent. Then the $i$th agent sends the encrypted messages $\mathcal{E}(A_{ui}x_i^k+r_ic)$ and $\mathcal{E}(A_{gi}x_i^k+s_id)$ to the SO, where $\mathcal{E}(\cdot)$ denotes any additively homomorphic encryption scheme, e.g., SingleMod encryption \cite{ryan2014doublemod}, Benaloh cryptosystem \cite{cohen1985robust}, and Paillier encryption \cite{paillier1999public}. Then the SO collects the encrypted messages from all the agents and performs arithmetic operations on the ciphertexts to obtain $\mathcal{E}(\sum_{i=1}^{n} A_{ui} x_i + c)$ and $\mathcal{E}(\sum_{i=1}^{n} A_{gi} x_i + d)$. This procedure will fully protect the private information $x_i$, because the SO is only dealing with ciphertext. As aforementioned,  $\sum_{i=1}^{n} A_{ui} x_i + c$ and $\sum_{i=1}^{n} A_{gi} x_i + d$ are indispensable components for the calculation of the subgradients in \eqref{6} and \eqref{7}, respectively, which would be used for primal and dual updates. Then the SO sends $\mathcal{E}(\sum_{i=1}^{n} A_{ui} x_i + c)$ and $\mathcal{E}(\sum_{i=1}^{n} A_{gi} x_i + d)$ to the agents. Finally, the agents decrypt the ciphertexts, obtain the plaintexts $\sum_{i=1}^{n} A_{ui} x_i + c$ and $\sum_{i=1}^{n} A_{gi} x_i + d$, and execute both primal and dual updates under plaintexts.



This novel paradigm is applicable to general primal-dual based algorithms. For example, the regularized primal-dual subgradient method (RPDS) \cite{koshal2011multiuser} takes the form of  
\begin{subequations}
\begin{align}
    x_{i}^{k+1} &=\Pi_{\mathbb{X}_{i}}(x_{i}^{k}-\alpha_{i} \nabla_{x_{i}} \mathcal{L}_{v, \epsilon}(x^{k}, \lambda^{k})),\\
    \lambda^{k+1} &= \Pi_{\mathbb{D}_{v}}(\lambda^{k}+\beta \nabla_{\lambda} \mathcal{L}_{v, \epsilon}(x^{k}, \lambda^{k})),
\end{align}
\end{subequations}
where $\alpha_i$ is the primal update step size, $\beta$ is the dual update step size, $\mathcal{L}_{v, \epsilon}(x, \lambda)=\mathcal{F}(x)+\frac{v}{2}\|x\|^{2}+\lambda^{\mathsf{T}} h(x)-\frac{\epsilon}{2}\|\lambda\|^{2}$ is the regularized Lagrangian function with regulation parameters $v$ and $\epsilon$, $k$ is the iteration index, and $\Pi_{\star}(\cdot)$ denotes the Euclidean projection. Though RPDS can guarantee the convergence, regularization errors inevitably exist. As an improvement, shrunken primal-dual subgradient (SPDS) developed in \cite{liu2013decentralized} can effectively eliminate regularization errors and speed up the convergence. The SPDS takes the form of
\begin{subequations}
\begin{align}
x_{i}^{k+1} &=\Pi_{\mathbb{X}_{i}}\left( \frac{1}{\tau_x} 
\Pi_{\mathbb{X}_{i}}\left( \tau_xx_{i}^{k}-\alpha_{i} \nabla_{x_{i}} \mathcal{L}\left(x^{k}, \lambda^{k}\right)\right)\right), \label{11a}\\
\lambda^{k+1} &= \Pi_{\mathbb{D}}\left( \frac{1}{\tau_{\lambda}}\Pi_{\mathbb{D}}\left(
\tau_{\lambda}\lambda^{k}+\beta \nabla_{\lambda} \mathcal{L}\left(x^{k}, \lambda^{k}\right)\right)\right),\label{11b}\
\end{align}
\end{subequations}
where $\tau_x$ and $\tau_{\lambda}$ denote the shrunken parameters for the primal and dual updates, respectively, and $\mathcal{L}(\cdot)$ denotes the relaxed Lagrangian defined in \eqref{5}. In what follows, we choose SPDS to integrate into the proposed privacy-preserving paradigm.

For simplicity of illustration, in this paper, we adopt the SingleMod cryptosystem which encrypts a plaintext integer message $z_r$ by
\begin{align}
    e = \mathcal{E}(z_r)=mw + z_r,
    \label{11}
\end{align}
where $e$ is the ciphertext, $m$ is a random integer, and $w$ is the private key which is a large prime number. Note that the coefficients and variables (e.g., $A_{ui}$ and $x_i$) in the decentralized calculations are real numbers, however $z_r$ has to be an integer in any cryptosystem. Therefore a real number $r$ is transformed into an integer $z_r$ by $z_r = 10^{\sigma} r$ where $\sigma$ denotes the preserved decimal fraction digits. A ciphertext can be transformed back to a real number by \cite{lu2018privacy}
\begin{align}
T_{\sigma, w}(z)= \begin{cases}
     z/ 10^{\sigma}, &\text { if } 0 \leq z \leq(w-1) / 2, \\ 
    (z-w) / 10^{\sigma}, &\text { if }(w+1) / 2 \leq z<w,
    \end{cases}
    \label{12}
\end{align}
where $z = e \bmod w$ and $\bmod$ denotes the modular operation. The SingleMod encryption is both additively homomorphic and multiplicatively homomorphic, indicating
\begin{subequations} \label{fully_homorphic}
\begin{align}
\mathcal{D}(\sum_{\ell=1}^{n} \mathcal{E}(z_{r\ell},w,m)) &=\sum_{\ell=1}^{n} z_{r\ell}, \label{13} \\
       \mathcal{D}(\prod_{\ell=1}^{n} \mathcal{E}(z_{r\ell},w,m))&=\prod_{\ell=1}^{n} z_{r\ell}, \label{14s} 
\end{align}
\end{subequations}
where $z_{r\ell}$ is the $\ell$th plaintext message and $\mathcal{D}(\cdot)$ denotes the decryption process. In \eqref{fully_homorphic}, summation and multiplication operations are performed on the ciphertexts, and when decrypted the output is the same as if those operations were performed on the plaintexts. In \cite{lu2018privacy}, the fully homomorphic properties of the chosen encryption scheme are mandated in designing a private key computing algorithm. In contrast, our proposed paradigm only requires additively homomorphic. The procedure of the proposed paradigm is presented via Algorithm \ref{alrogithm_1}.

\noindent {\bf{Remark 1:}} The proposed paradigm is still efficacious in special cases when $c=0$ or $d=0$. Take $c=0$ for example, the SO could select any random $c$ with an appropriate dimension, replace \eqref{8a} with a set of random numbers satisfying $r_1+\ldots+r_n = 0$, and send $r_ic$ to agent $i$, $\forall i=1,\ldots,n$.  \hfill $\blacksquare$

\noindent {\bf{Remark 2:}} Since all agents share the same private key, the privacy could be compromised if agent $i$ wiretaps the ciphertexts sent from agent $j$. To overcome this issue, a mandatory assumption was made in \cite{lu2018privacy} requiring that no agent can eavesdrop the communication links between other participants. Algorithm \ref{alrogithm_1} lifts this assumption by introducing randomized $r_ic$ and $s_id$ into the exchanged messages so that an agent cannot extract the true decision variables of other agents even if it eavesdrops the communication. \hfill $\blacksquare$

\noindent {\bf{Remark 3:}} Algorithm \ref{alrogithm_1} is applicable to both private key encryption scheme (e.g., SingleMod encryption scheme \cite{ryan2014doublemod}) and 
public key encryption scheme (e.g., Benaloh cryptosystem \cite{cohen1985robust} and  Paillier encryption scheme \cite{paillier1999public}). \hfill $\blacksquare$

\begin{algorithm} 
\caption{Cryptography-based privacy-preserving decentralized optimization paradigm}
\begin{algorithmic}[1]
\State All agents agree on a private key $w$.
\State Agents initialize primal and dual variables, tolerance $\epsilon_0$, iteration counter $k{=}0$, and maximum iteration $k_{max}$.

\While{ \textit{ $\epsilon > \epsilon_0$ and $k < k_{max}$}}

\State SO generates a set of parameters satisfying $\sum_{i=1}^n r_i^k = 1$ and 
$\sum_{i=1}^n s_i^k {=} 1$, then sends $r_i^kc$ and $s_i^kd$ to the $i$th agent.

\State  All agents adopt SingleMod encryption for $ A_{ui} x_i^k + r_i^kc$ and $ A_{gi} x_i^k + s_i^kd$ using \eqref{11}, then send $\mathcal{E}(A_{ui} x_i^k + r_i^kc)$ and $\mathcal{E}(A_{gi} x_i^k + s_i^kd)$  to the SO.

\State  SO firstly collects the encrypted messages from all agents, then sums the received messages over ciphertext using \eqref{13} to obtain $\mathcal{E}(\sum_{i=1}^{n} A_{ui} x_i^k + c)$ and $\mathcal{E}(\sum_{i=1}^{n} A_{gi} x_i^k + d)$. Then the SO sends  $\mathcal{E}(\sum_{i=1}^{n} A_{ui} x_i^k + c)$ and $\mathcal{E}(\sum_{i=1}^{n} A_{gi} x_i^k + d)$ to the $i$th agent.

\State  Agent $i$ receives and decrypts  $\mathcal{E}(\sum_{i=1}^{n} A_{ui} x_i^k + c)$ and $\mathcal{E}(\sum_{i=1}^{n} A_{gi} x_i^k + d)$ from the SO with the private key  $w$, then converts the integer to real number by using \eqref{12}.

\State  Each agent $i$ updates the primal variable $x_i$ by \eqref{11a} and the dual variable $\lambda$ by \eqref{11b}
over plaintext.

\State Calculate the error $\epsilon$. 

\State  $k=k+1$.
\EndWhile
\end{algorithmic}
\label{alrogithm_1}
\end{algorithm}

\subsection{Privacy Analysis}
Due to the limit of the space, in this section, we concisely present high-level analyses of the privacy preservation advantages of the proposed paradigm. Systematic analysis and theoretical proofs will be provided in the journal extension of this work. Particularly, three types of adversaries are considered: 
\begin{itemize}
  \item  \emph{Honest-but-curious agents} -- an honest-but-curious agent is only interested in solving the optimization problem, but does not tamper with the algorithm. However, it may observe the intermediate or input/output data to infer the private information of other participating agents;
  \item \emph{External eavesdroppers} -- external eavesdroppers launch attacks by wiretapping and intercepting exchanged messages between agents and the SO;
  \item \emph{SO} -- SO may attempt to infer the decision variables of the agents by collecting the received data from the agents.
\end{itemize}

 First, by introducing the random parameters $r_i$ and $s_i$ generated by the SO: (1) The public coefficients $c$ and $d$ are kept private to the SO, because only the randomized $r_ic$ and $s_id$ are transmitted to the agents; (2) The $i$th agent passes $\mathcal{E}(A_{ui}x_i^k+r_ic)$ and $\mathcal{E}(A_{gi}x_i^k+s_id)$ to the SO. When \textit{honest-but-curious} agents wiretap those messages, the only accessible information (from direct decryption) is $A_{ui}x_i^k+r_ic$ and $A_{gi}x_i^k+s_id$, which are randomized by $r_ic$ and $s_id$. Comparing with the private key encryption scheme in \cite{lu2018privacy}, Algorithm \ref{alrogithm_1} achieves enhanced security by protecting the privacy of an agent from other \textit{honest-but-curious} agents. 

Algorithm \ref{alrogithm_1} also effectively blocks information leakage to \textit{external eavesdroppers}. This is because all messages transmitted within the proposed framework are either encrypted or randomized. Without the secret private key, external eavesdroppers are not able to decrypt those messages. Moreover, even when the private key is accidentally available to the external eavesdroppers, owing to the extra randomization, the true message still cannot be revealed. This security always holds unless in an extreme case where the \textit{external eavesdroppers} wiretap the received data of all the agents. 

Finally, though the SO could be a trustworthy participant, the agents may be unwilling to share their private decision variables to the SO. In the worst-case scenario, SO may collect the agents' data sequence to infer their true states and intentions. In \cite{lu2018privacy}, due to the lack of a fully homomorphic encryption scheme for integers that is both efficiently implementable and semantically secure, the agents' privacy could be compromised if the SO fully knows the update rule of $x$ and uses the collected data sequence to estimate the agents' decision variables. This issue is resolved in Algorithm \ref{alrogithm_1} as the SO has no knowledge about the update rule of the agents.



\section{Simulation Results}
In this section, we verify the correctness and efficiency of the proposed privacy-preserving paradigm with a numerical example and a traffic congestion optimization problem. In both examples, the encryption error is defined by 
\begin{equation}
    P_e^k = \sum_{i=1}^n \| \hat{x}_i^k - x_i^k\|_2,
\end{equation}
where $\hat{x}_i^k$ is the solution of Algorithm \ref{alrogithm_1} and $x_i^k$ is obtained by solving the optimization problem without privacy preservation. The optimality gap at the $k$th iteration is defined as
\begin{equation}
    G_e^k = \sum_{i=1}^n \| \hat{x}_i^k - x_i^*\|_2,
\end{equation}
where $x_i^*$ is the optimizer and the convergence tolerance is set to be $\epsilon_0 {=} 10^{-4}$. The shrunken parameters were empirically chosen as $\tau_x {=} \tau_{\lambda} {=}0.98$ and the precision level was set as $\sigma {=} 3$. The private key is normally extremely large in real applications for security purpose, e.g., in the magnitude of $2^{2000}$, however, for clear presentation, we chose the magnitude of $2^{50}$. 

\subsection{A Numerical Example}
Consider a case with two agents and one SO, and the optimization problem can be written as  
\begin{equation}
\begin{aligned}
& \underset{x}{\text{min}} & & {\mathcal{F}(x)} \\
& \: \text{s.t.} & & x_{i} \in [0,1], \quad  i=1,2,\\
& & &  A_{g1} x_1 + A_{g2} x_2 + d \leq 0,
\end{aligned}
\label{14}
\end{equation}
where ${\mathcal{F}(x)}$ denotes the objective function which is given by 
\begin{align}
    {\mathcal{F}(x)} = &\frac{1}{2}\| A_{u1} x_1 + A_{u2} x_2 +c \|_2^2
    {+} (A_{q1}x_1)^{\mathsf{T}}(A_{q1}x_1) {+} A_{l1} x_1 \nonumber\\&+ C_{t1}+ (A_{q2} x_2)^{\mathsf{T}}(A_{q2} x_2)
     +A_{l2} x_2 +C_{t2}.
\end{align}
The coefficients were chosen as
$A_{u1} = [\begin{smallmatrix} -1&0\\ 1&-0.5 \end{smallmatrix}]$, $A_{u2} = [\begin{smallmatrix} 0&-2\\ 0&-10 \end{smallmatrix}]$,
$c {=} [\begin{smallmatrix} 1\\ 1 \end{smallmatrix}]$,
$A_{q1} {=} [\begin{smallmatrix} 1&0\\ 1&1 \end{smallmatrix}]$, $A_{q2} {=} [\begin{smallmatrix} 0&1\\ 1&1 \end{smallmatrix}]$, $A_{l1} {=} [\begin{smallmatrix} 1&1 \end{smallmatrix}]$, $A_{l2} {=} [\begin{smallmatrix} 1&0 \end{smallmatrix}]$, $C_{t1} {=}1$,  $C_{t2} {=}0$, $A_{g1} = [\begin{smallmatrix} 1&0\\ 1&-1 \end{smallmatrix}]$, $A_{g2} = [\begin{smallmatrix} 0&1\\ -1&-1 \end{smallmatrix}]$, and $d = [\begin{smallmatrix} -1\\ 1 \end{smallmatrix}] $. Note that the coefficients were carefully chosen to be representative, i.e., signed real numbers. The SO is trustworthy and has access to the global coefficients $A_{u1}$, $A_{u2}$ ,$A_{g1}$, $A_{g2}$, $c$, and $d$. Each agent only has the global coefficients that are related to itself and the individual coefficients included in $\mathcal{I}_{fi}$, i.e., agent 1 has access to $A_{u1}$, $A_{g1}$, and $\mathcal{I}_{f1} = \{A_{q1},A_{l1},C_{t1}\}$. The coefficients assignment is given by Table \ref{table_numerical_coe_assign}.
\vspace*{-2mm}
\begin{table}[!htb]
\caption{Coefficients assignment of the numerical example}
\vspace*{-3mm}
\label{table_numerical_coe_assign}
\begin{center}
\begin{tabular}{c|l}
\hline
Participant name & Coefficients held \\
\hline
SO & $A_{u1},A_{u2},A_{g1},A_{g2},c,d$ \\
 Agent 1 & $\mathcal{I}_{f1},A_{u1},A_{g1}$    \\
 Agent 2 & $\mathcal{I}_{f2},A_{u2},A_{g2}$   \\
\hline
\end{tabular}
\end{center}
\vspace*{-2mm}
\end{table}

\begin{figure}[!htb]
\vspace*{-2mm}
    \centering
    \includegraphics[width=0.34\textwidth, trim = 26mm 40mm 26mm 48mm, clip]{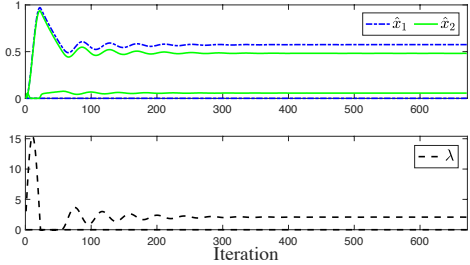}
    \vspace*{-3mm}
    \caption{Convergence of the primal and dual variables.}
    \label{numerical_primal_dual_convergence}
\end{figure}

\begin{figure}[!htbp]
\vspace*{-2mm}
    \centering
    \includegraphics[width=0.34\textwidth,trim =26mm 40mm 26mm 48mm, clip]{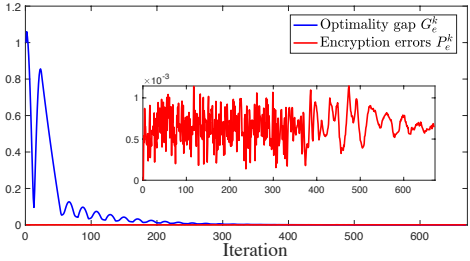}
    \vspace*{-3mm}
    \caption{Encryption error and optimality gap (precision level $\sigma = 3$).}
    \label{numerical_convergence_error}
\end{figure}

\begin{figure}[!htbp]
    \centering
    \includegraphics[width=0.35\textwidth, trim = 26mm 35mm 26mm 43mm, clip]{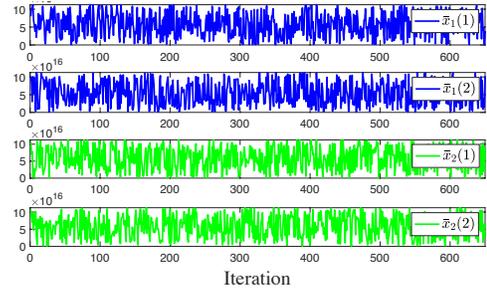}
    \vspace*{-3mm}
    \caption{Secure information exchange between agents and the SO.}
    \label{secure_information}
\end{figure}

The true optimizers of \eqref{14} are $x_1^* = [0, 0.5750]^{\mathsf{T}}$ and
$x_2^* = [0.4814, 0.0564]^{\mathsf{T}}$. We then solve \eqref{14} using Algorithm \ref{alrogithm_1}. The primal step sizes were uniformly chosen as $\alpha = 5 \times 10^{-3}$ and the dual step size was set to $\beta = 2$. Fig. \ref{numerical_primal_dual_convergence}  presents that the primal and dual variables converged in about 600 iterations.   Fig. \ref{numerical_convergence_error} shows the optimality gap of the primal variables and the encryption errors caused by real number and integer transformation.  The encryption errors are strictly less than 0.01 and therefore 
comply with the precision level $\sigma = 3$. Fig. \ref{secure_information} shows the exchanged messages $\bar{x}_i = \mathcal{E}(A_{ui} x_i + r_ic)$ between agents and the SO.  Therefore, the encrypted and randomized messages in Fig. \ref{secure_information} are useless for any adversaries to extract any private information of the agents.

\subsection{Traffic Congestion Optimization}
\begin{figure}[!htbp]
\vspace*{-2mm}
    \centering
    \includegraphics[width=0.2\textwidth,trim = 5mm 5mm 5mm 10mm, clip]{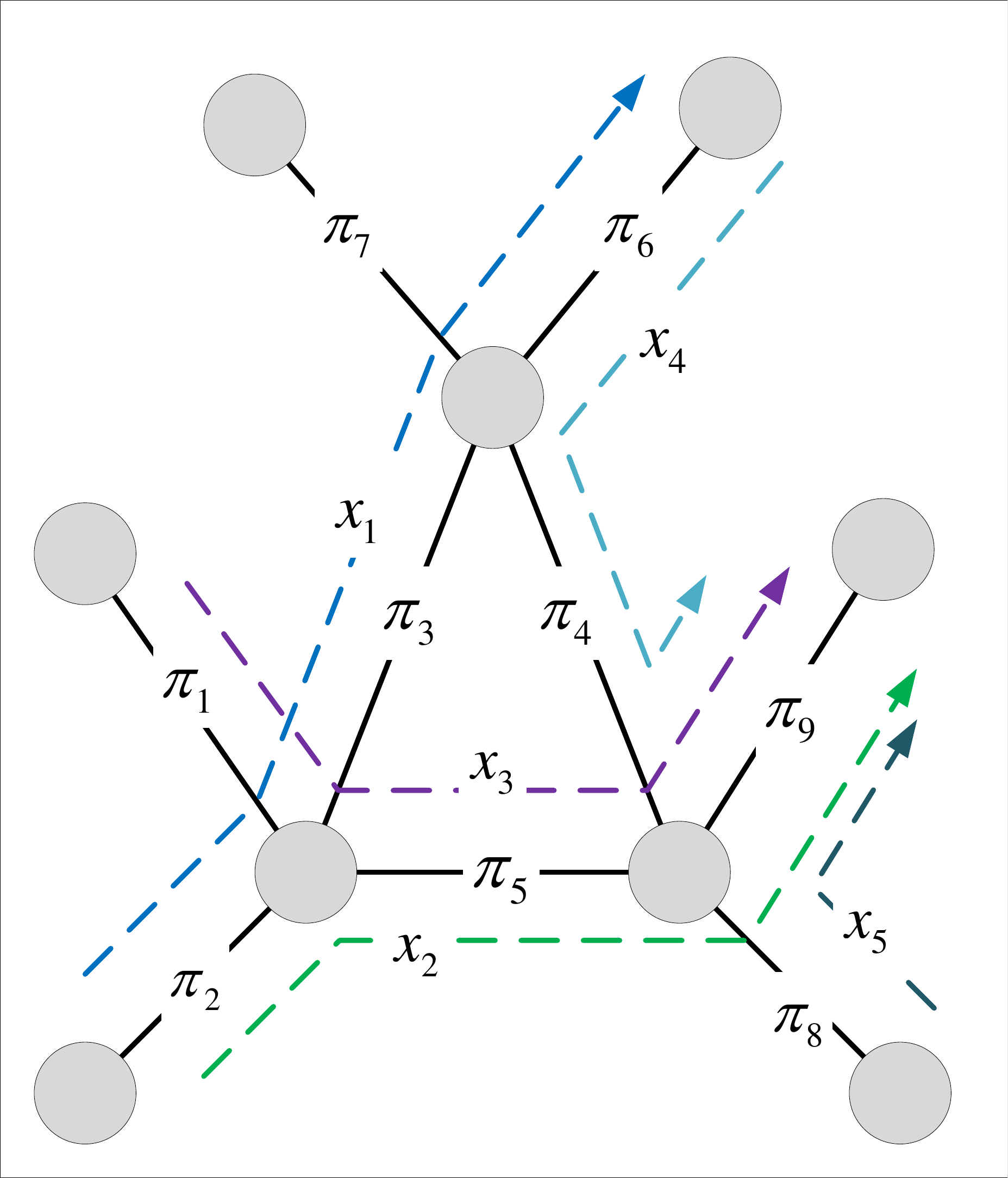}
    \vspace*{-3mm}
    \caption{A transportation network with 9 links and 5 agents \cite{koshal2011multiuser}.}
    \label{traffic_flow}
\end{figure}

The traffic congestion problem was formulated over a transportation network with $\mathcal{N}$ agents and $\mathcal{L}$ links.
Suppose each agent $i$ travels along a route with the transmission rate $x_i$, then the shared paths between all agents arise traffic congestion. An example of a transportation network with $\mathcal{N}=5$ users and $\mathcal{L}=9$ links is shown in Fig. \ref{traffic_flow}.  This example was used in \cite{koshal2011multiuser} to demonstrate the effectiveness of RPDS. Let $\pi_j$ be the $j$th link, $x = [x_1 \cdots x_\mathcal{N}]^{\mathsf{T}}$, and $A \in \mathbb{R}^{\mathcal{L}\times \mathcal{N}}$ denote the link-route incidence matrix, i.e., $A_{ji}=1$ if link $j$ is on the path of agent $i$, and $A_{ji}=0$ otherwise. Consequently, the transportation congestion optimization problem can be formulated as 
\begin{subequations} \label{eq:15}
\begin{align}
    &\underset{x}{\text{min}} \quad { c(x) + \sum_{i=1}^{\mathcal{N}} f_i(x_i)} \label{eq:12sa}\\
& \: \text{s.t.}  \ \quad x \geq 0, \label{eq:12sb}\\
& \,\,\, \qquad \sum_{i=1}^{\mathcal{N}} A^{col}_ix_i -b \leq 0. \label{eq:12sc}
\end{align}
\end{subequations}
The coupled congestion cost in \eqref{eq:12sa} is defined as 
\begin{equation}
    c(x) = \sum_{i=1}^{\mathcal{N}} \sum_{l \in \mathcal{L}}x_{li}\sum_{q=1}^{\mathcal{N}}x_{lq} = \| \sum_{i=1}^{\mathcal{N}} A^{col}_ix_i\|_2^2,
    \label{33ss}
\end{equation}
where $x_{li}$ is the flow of agent $i$ on link $l$ and $A^{col}_i$ denotes the $i$th column of $A$. The local cost of agent $i$ is given by \cite{kelly1998rate}
\begin{equation}
    f_i(x_i) = -k_i\log (1+x_i).
    \label{32sn}
\end{equation}
The global constraint \eqref{eq:12sc} limits the agents' traffic rates over the network where $b$ is the link capacity vector with $b_l$ denoting the maximum traffic capacity of link $l$. Table \ref{table_example_congestion_flow} concludes the traffic flow in the network and $k_i$ in the cost function. 
\begin{table}[!htb]
\caption{Traffic Network and Agent Data}
\vspace*{-3mm}
\label{table_example_congestion_flow}
\begin{center}
\begin{tabular}{c|l|l}
\hline
Agent name & Links traversed & $k_i$\\
\hline
 1 & $\lambda_2,\lambda_3,\lambda_6$    & 10 \\
 2 & $\lambda_2,\lambda_5,\lambda_9$   & 0 \\
 3 & $\lambda_1,\lambda_5,\lambda_9$   & 10\\
 4 & $\lambda_6,\lambda_4,\lambda_9$  & 10 \\
 5 & $\lambda_8,\lambda_9$  & 10 \\
\hline
\end{tabular}
\end{center}
\vspace*{-5mm}
\end{table}

\begin{figure}[!htbp]
\vspace*{-3mm}
    \centering
    \includegraphics[width=0.35\textwidth,trim = 26mm 40mm 26mm 47mm, clip]{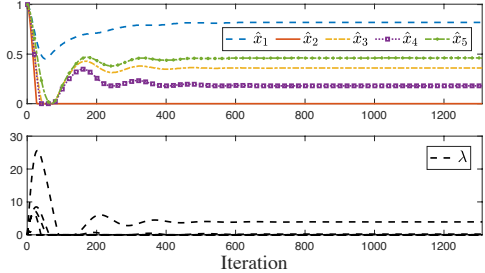}
    \vspace*{-3mm}
    \caption{Convergence of the primal and dual variables.}
    \label{congestion_primal_dual}
\end{figure}

\begin{figure}[!htbp]
    \centering
    \includegraphics[width=0.34\textwidth, trim = 26mm 40mm 26mm 47mm, clip]{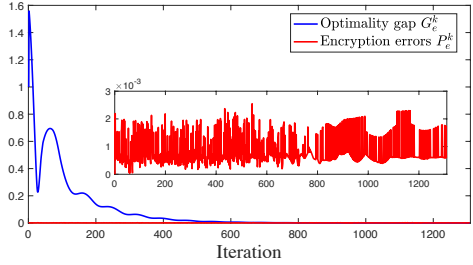}
    \vspace*{-3mm}
    \caption{Encryption error and optimality gap (precision level $\sigma = 3$).}
    \label{congestion_convergence_errors}
\end{figure}

The link capacity vector was set as $b = \bm{1}_9$. The primal step sizes were uniformly chosen as $\alpha = 10^{-3}$ and the dual step size was set to $\beta = 0.5$. Fig. \ref{congestion_primal_dual} shows that the primal and dual variables converged in about 1000 iterations. Fig. \ref{congestion_convergence_errors} shows that the proposed paradigm has converged to the optimal solutions and the encryption errors are strictly smaller than 0.01 during all the iterations.  Fig. \ref{congestion_primal_dual} and Fig. \ref{congestion_convergence_errors} verify that the proposed paradigm well inherits the feasibility and optimality of SPDS. 
\vspace*{-6mm}
\vspace*{-2mm}



\section{Conclusions}
This paper developed a novel privacy-preserving decentralized paradigm for a family of strongly coupled cooperative optimization problems. The proposed paradigm is scalable to the agent population size and achieves privacy preservation of all the participants. The privacy of both the agents and the SO was analyzed, and the privacy preservation is guaranteed under a range of adversaries. The efficiency and efficacy of the proposed paradigm were verified via a numerical example and a traffic congestion optimization problem. The future research directions of this paper would include (1) considering the case where the SO is not a trustworthy party and does not possess any knowledge of the optimization problem, (2) theoretically analyzing the security of the proposed paradigm, and (3) verifying the industrial values of the proposed paradigm. In specific, an experimental platform with a cluster of micro controllers could be designed to demonstrate ICPS applications.

\bibliographystyle{IEEEtran}
\bibliography{ICPS}


\end{document}